\def\mycmd{0} 
\if\mycmd0
\documentclass{article} 
\else
\documentclass{elsarticle} 
\fi

\usepackage{amsmath,amssymb}
\usepackage{mathtools} 
\usepackage{fullpage}  
\usepackage{todonotes} 
\usepackage{subcaption} 
\usepackage{booktabs}  
\usepackage{tikz}
\usepackage{pgfplots}
\usepackage{cancel} 
\pgfplotsset{compat=1.15}
\usepgfplotslibrary{colorbrewer}
\usepackage{csvsimple}
\usepackage{float}  
\makeatletter
\renewcommand{\fps@figure}{H}
\makeatother

\newtheorem{theorem}{Theorem}

\newcommand{\displacement}{\mathbf{d}}
\newcommand{\mat}[1]{\mathbf{#1}}
\newcommand{\ten}[1]{\mathbf{#1}}
\renewcommand{\vec}[1]{\mathbf{#1}}
\DeclareMathOperator{\dive}{div}
\DeclareMathOperator{\tr}{tr}
\DeclareMathOperator{\cof}{cof}
\newcommand{\grad}{\nabla}

\newcommand{\newton}{\textbf{NK}}
\newcommand{\inexactnewton}{\textbf{iNK}}
\newcommand{\bfgs}{\textbf{B}}
\newcommand{\inexactbfgs}{\textbf{iB}}
\newcommand{\Pone}{$\mathbb P_1$\,}
\newcommand{\Ptwo}{$\mathbb P_2$\,}

\newcommand{\newtoncolor}{blue!70!white}
\newcommand{\inexactnewtoncolor}{red!80!black}
\newcommand{\bfgscolor}{yellow!50!green!80!black}
\newcommand{\inexactbfgscolor}{yellow!70!black}

\pgfplotsset{every axis plot/.append style={line width=1.1pt}}
\newenvironment{plot}[3][]
      {\begin{tikzpicture}
      \begin{axis}[width=\textwidth, xlabel=#2, ylabel=#3, tick label style={font=\footnotesize}, legend style={font=\tiny, draw=none, fill opacity=0.5, text opacity = 1,row sep=2pt},legend pos=north west, legend columns=1, #1]
      }
      {%
      \end{axis}
      \end{tikzpicture}
      }

\newenvironment{semilogxplot}[3][]
      {\begin{tikzpicture}
      \begin{semilogxaxis}[width=\textwidth, xlabel=#2, ylabel=#3, tick label style={font=\footnotesize}, legend style={font=\tiny, draw=none, fill opacity=0.5, text opacity = 1,row sep=2pt},legend columns=1,legend pos=north west, #1]
      }
      {%
    \end{semilogxaxis}
      \end{tikzpicture}
      }
\newenvironment{semilogyplot}[3][]
      {\begin{tikzpicture}
      \begin{semilogyaxis}[width=\textwidth, xlabel=#2, ylabel=#3, tick label style={font=\footnotesize}, legend style={font=\tiny, draw=none, fill opacity=0.5, text opacity = 1,row sep=2pt},legend columns=1,legend pos=north west, #1]
      }
      {%
    \end{semilogyaxis}
      \end{tikzpicture}
      }

\newenvironment{loglogplot}[3][]
      {\begin{tikzpicture}
      \begin{loglogaxis}[width=\textwidth, xlabel=#2, ylabel=#3, tick label style={font=\footnotesize}, legend style={font=\tiny, draw=none, fill opacity=0.5, text opacity = 1,row sep=2pt},legend columns=1, #1]
      }
      {%
    \end{loglogaxis}
      \end{tikzpicture}
      }

\usepackage{listings}
\usepackage{xcolor}
\definecolor{codegreen}{rgb}{0,0.6,0}
\definecolor{codegray}{rgb}{0.5,0.5,0.5}
\definecolor{codepurple}{rgb}{0.58,0,0.82}
\definecolor{backcolour}{rgb}{0.95,0.95,0.92}
\lstdefinestyle{mystyle}{
  backgroundcolor=\color{backcolour}, commentstyle=\color{codegreen},
  keywordstyle=\color{magenta},
  numberstyle=\tiny\color{codegray},
  stringstyle=\color{codepurple},
  basicstyle=\ttfamily\footnotesize,
  breakatwhitespace=false,         
  captionpos=b,                    
  keepspaces=true,                 
  numbers=none,                    
  numbersep=5pt,                  
  showspaces=false,                
  showstringspaces=false,
  showtabs=false,                  
  tabsize=2,
  framerule=1.5pt,
  rulecolor=\color{red!60!black}
}

\begin{document}
\if\mycmd0
\title{Parallel inexact Newton-Krylov and quasi-Newton solvers for nonlinear elasticity}
\author{Nicolás A. Barnafi\thanks{Department of Mathematics, Universit\`a di Pavia, Via Ferrata 1, 27100 Pavia, Italy. \texttt{\{nicolas.barnafi, luca.pavarino\}@unipv.it}} \and Luca F. Pavarino\footnotemark[1] \and Simone Scacchi\thanks{Department of Mathematics, Universit\`a di Milano,
Via Saldini 50, 20133 Milano, Italy. \texttt{simone.scacchi@unimi.it}}}
\date{}
\maketitle
\begin{abstract}
	In this work, we address the implementation and performance of inexact Newton-Krylov and quasi-Newton algorithms, more specifically the BFGS method, for the solution of the nonlinear elasticity equations, and compare them to a standard Newton-Krylov method. This is done through a systematic analysis of the performance of the solvers with respect to the problem size, the magnitude of the data and the number of processors in both almost incompressible and incompressible mechanics. We consider three test cases: Cook's membrane (static, almost incompressible), a twist test (static, incompressible) and a cardiac model (complex material, time dependent, almost incompressible). Our results suggest that quasi-Newton methods should be preferred for compressible mechanics, whereas inexact Newton-Krylov methods should be preferred for incompressible problems. We show that these claims are also backed up by the convergence analysis of the methods. In any case, all methods present adequate performance, and provide a significant speed-up over the standard Newton-Krylov method, with a CPU time reduction exceeding 50\% in the best cases.
\end{abstract}
\else
\begin{frontmatter}

 \title{Parallel inexact Newton-Krylov and quasi-Newton solvers for
nonlinear elasticity}

 \author[1]{Nicolás Barnafi}
 \author[2]{Simone Scacchi}
 \author[1]{Luca Pavarino}

\address[1]{Department of Mathematics, Universit\`a di Pavia, Via Ferrata 1, 27100 Pavia, Italy. \texttt{\{nicolas.barnafi, luca.pavarino\}@unipv.it}}
\address[2]{Department of Mathematics, Universit\`a di Milano,
Via Saldini 50, 20133 Milano, Italy. \texttt{simone.scacchi@unimi.it}}

\date{\today}
 
\begin{abstract}
In this work, we address the implementation and performance of inexact Newton-Krylov and quasi-Newton algorithms, more specifically the BFGS method, for the solution of the nonlinear elasticity equations, and compare them to a standard Newton-Krylov method. This is done through a systematic analysis of the performance of the solvers with respect to the problem size, the magnitude of the data and the number of processors in both almost incompressible and incompressible mechanics. We consider three test cases: Cook's membrane (static, almost incompressible), a twist test (static, incompressible) and a cardiac model (complex material, time dependent, almost incompressible). Our results suggest that quasi-Newton methods should be preferred for compressible mechanics, whereas inexact Newton-Krylov methods should be preferred for incompressible problems. We show that these claims are also backed up by the convergence analysis of the methods. In any case, all methods present adequate performance, and provide a significant speed-up over the standard Newton-Krylov method, with a CPU time reduction exceeding 50\% in the best cases.
\end{abstract}


\end{frontmatter}
\fi

\section{Introduction}

Nonlinear elasticity is a continuum framework for modeling the deformation of an elastic body, and has a large spectrum of applications, ranging from industrial to academic. The flexibility of continuum mechanics resides in its abstract formulation starting from fundamental principles, which allows for an accurate representation of many different materials through constitutive modeling and complex boundary conditions \cite{Holzapfel2002489}. After a precise physical phenomenon has been devised and modeled, the resulting equations yield a complex system of nonlinear Partial Differential Equations (PDEs). Limited progress has been achieved in terms of its numerical analysis \cite{carstensen2004priori,auricchio2013approximation}, so it is difficult to know a-priori efficient strategies for its numerical approximation. Indeed, the problem is not convex but only polyconvex \cite{ciarlet2021mathematical}, which deteriorates the performance of many well established numerical methods. The gold standard is Newton's method \cite{zienkiewicz1977finite}, consisting in a second order approximation of the associated variational principle (or a first order approximation of the Euler-Lagrange equations, referred to as the Newton-Raphson method). This method is popular mainly due to its robustness and low iteration count, as it converges quadratically whenever a good initial guess of the solution is considered \cite{wright1999numerical}. Its main drawback is that it requires the repeated assembly of the Jacobian matrix, associated to the tangent problem, which usually becomes a bottleneck of the solution process (see \cite{jansson2011adaptive} for an example in CFD). This can be more evidently appreciated when using higher order approximations, required for avoiding numerical locking effects \cite{ern2013theory}.

In general, not much attention has been given in the community to other methods, such as linearly convergent (gradient descent \cite{ruder2016overview}, Richardson \cite{saad2003iterative}) and superlinearly convergent (quasi-Newton \cite{wright1999numerical}, inexact Newton \cite{dembo1982inexact}) ones. These methods yield a higher iteration count due to their reduced convergence rate, but can potentially present a drastically overall reduced computational complexity. Indeed, the use of different nonlinear solvers has shown enormous speed-ups in other problems \cite{smith2006intelligent,borregales2017robust}, where the properties of the equations have been exploited to devise a more adequate nonlinear solution method. Some recent works have numerically explored quasi-Newton methods for the mechanics \cite{linge2005solving,gelin1988use,liu2017quasi}, with their main focus being on the nonlinear iterations and CPU time. This type of study does not take into account the linear system involved in each iteration of the method, and thus is not conclusive with respect to the applicability of these methods in an HPC setting, where problems with millions of degrees of freedom need to be solved, possibly for many time steps. A numerical method can be deemed adequate in such case if at least is satisfies the following requirements: i) the nonlinear and linear iterations incurred during the solution procedure are independent of the mesh size used, ii) the method is robust with respect to the model parameters, at least in the scenarios of interest and iii) the method is (strongly) scalable, meaning that the overall solution time improves when more processors are used. We finally highlight \cite{weiser2007affine}, where a novel conjugate Newton method was proposed and analysed for nonlinear elasticity.

The scope of this work is to provide a first step in the adequate usage of alternatives to Newton's method, where we focus in superlinearly convergent methods as they provide an excellent overall computational complexity \cite{wright1999numerical}. We stress that all gradient descent algorithms we tested on a preliminary phase failed, and this is indeed consistent with the literature: there are no works, up to our knowledge, of descent nor fixed point algorithms for nonlinear mechanics. We focus on three tests: i) Cook's membrane test, an almost incompressible problem, ii) a twist problem, where we test an incompressible material and iii) a cardiac modeling test, where we use an idealized left ventricle geometry to model a human heartbeat. The work is structured as follows: in Section \ref{section:hyperelasticity} we review the mechanics problem and fix some notation, in Section \ref{section:solvers} we describe the nonlinear solvers to be used throughout this study, in Section \ref{section:tests} we define the tests in which the methods will be tested together with a thorough description of the linear solver involved, in Section \ref{section:theory} we recall the convergence results of the methods and discuss what their expected performance should be, in Section \ref{section:results} we show and comment the results obtained and we conclude with a discussion of the results in Section \ref{section:discussion}.

\section{The hyperelasticity problem}\label{section:hyperelasticity}
Consider a connected domain $\Omega\subset \mathbb R^{d\in\{2,3\}}$ that represents the geometry to be deformed, and define its Dirichlet and Neumann boundaries as $\partial \Omega_D$ and $\partial \Omega_N$ respectively, such that $\overline{\partial\Omega} = \overline{\partial\Omega_D}\cup \overline{\partial\Omega_N}$. We look for a displacement $\displacement: [0,T]\times \Omega \to \mathbb R^d$ such that it solves the momentum conservation equation
	\begin{equation}\label{eq:momentum strong}
		\begin{aligned}
			\ddot \displacement - \dive \ten P &= \vec f &&\text{in $(0,T)\times\Omega$},\\
			\displacement &= \displacement_0 &&\text{in $\{0\}\times\Omega$}, \\
			\dot{\displacement} &= \vec v_0 &&\text{in $\{0\}\times \Omega$}, \\
			\ten P\vec \vec n &= \vec t &&\text{on $[0,T)\times\partial\Omega_N$}, \\
			\displacement &= \displacement_D &&\text{on $[0,T)\times\partial\Omega_D$},\\
		\end{aligned}
	\end{equation}
where $\displacement_0$ is the initial displacement, $\vec v_0$ is the initial velocity,$\vec t$ is the surface traction, $\displacement_D$ is a distributed load, $\dot{()}$ denotes a time derivative, and $\ten P$ is known as the Piola stress tensor, usually obtained from a predefined Helmholtz potential $\Psi$ such that
	$$ \ten P(\ten F) = \frac{\partial \Psi}{\partial \ten F}(\ten F), $$
where $\ten F = \ten I + \grad \displacement$. The existence of such a potential is actually a hypothesis, and whenever it exists we refer to the material as a \emph{hyperelastic material}.
Approximating this problem by neglecting the inertial term $\ddot{\displacement}$ results in the well-known \emph{static} mechanics, or quasi-static whenever one of the problem's data is time-dependent. If we require the volume to be conserved through $J \coloneqq \det \ten F =1$, we obtain an \emph{incompressible} elasticity problem. This can be instead approximated by further penalizing the deviation of $J-1$ from 0, which results in a modified potential
	$$\Psi(\ten F) = \Psi_\text{sol}(\bar{\ten F}) + \Psi_\text{vol}(J),$$
where the isochoric component $\bar{\ten F}=J^{-1/d} \ten F$ is such that $\det \bar{\ten F} = 1$. This separates the energetic contribution of volumetric deformation, and is referred to as almost (or quasi-) incompressibility whenever the term $\Psi_\text{vol}$ yields $J\approx 1$. More details 
can be found in \cite{Holzapfel2002489}.

\section{Numerical solvers}\label{section:solvers}
In this section we briefly review the nonlinear methods we consider, namely Newton-Krylov and quasi-Newton methods. Details regarding their geometric motivation can be found in \cite{wright1999numerical}, whereas the convergence properties are detailed in Section \ref{section:theory}. Our study is mainly motivated by the use of these methods in an HPC infrastructure, so the use of iterative solvers within the nonlinear method is not only inevitable, but also desirable due to their well-established parallel performance. Inspired by this, we have classified our methods according to the degree of exactness of the linear solver it uses, so that a fully inexact method will simply consider the action of the associated preconditioner instead of solving a linear system. Instead, a quasi-exact\footnote{We emphasize that all methods are inexact given their iterative nature. Still, a sufficiently low tolerance can be regarded as a sufficiently exact solver, so we refer to this as quasi-exact.} method will consider an accurate linear solver given by small tolerances, i.e. $10^{-10}$ and $10^{-8}$ absolute and relative tolerances, respectively. Now we provide more details regarding the methods, and in particular emphasize where the linear solver is used.

\begin{description}
	\item[Newton-Krylov methods.] These methods can be seen as either a minimization procedure or as a root finding algorithm, the latter often referred to as Newton-Raphson method. We present it as a root finding method. Consider Equation \eqref{eq:momentum strong} in quasi-static, residual form
		$$ \vec R(\displacement) \coloneqq \vec f -\dive \ten P(\ten F) = \vec 0, $$
		with $\ten P = \frac{\partial \Psi}{\partial \ten F}$ for a given Helmholtz potential $\Psi$, and consider an initial iterate $\displacement^0$. Then, given a previous iteration $\displacement^{n-1}$, the next one is given as the solution $\delta \displacement^n$ of the linearized problem
		\begin{equation}\label{eq:newton}
			\partial_{\displacement}\vec R(\displacement^{n-1})[\delta \displacement^n] = -\vec R(\displacement^{n-1}),
		\end{equation}
		where $\partial_{\displacement}$ stands for the Frechét derivative with respect to $\displacement$, and the update is then given by $\displacement^n = \displacement^{n-1} + \delta\displacement^n$. Equation \eqref{eq:newton} gives rise to a linear system of equations often referred to as the tangent problem or simply the linearized problem, which is solved by means of an iterative Krylov space method, accelerated by a preconditioner, see e.g. \cite{Cai1994,Cai1998}. In practice, the most used one is referred to as Newton-MG, meaning that problem \eqref{eq:newton} is solved with an iterative method, preconditioned by an algebraic multigrid preconditioner (AMG) \cite{hackbusch2013multi,colli.2015}.
	    
	    \item[Inexact Newton-Krylov methods.] Problem \eqref{eq:newton} is an approximation of the actual equation $\vec R(\displacement)=0$, and thus it might not be true that an accurate solution of the tangent problem will also yield an accurate solution of the original equation. In fact, this usually gives rise to over-solving the problem \cite{eisenstat1994globally} in the first iterations, which motivates the use of inexact solvers for \eqref{eq:newton}, i.e. to use large tolerances for the solution of the tangent problem. In addition, the quality of the linearization improves as the iterates are closer to the solution, which motivates the use of adaptive (relative) tolerances. In particular, we use the Einstat-Walker strategy \cite{eisenstat1996choosing}, given by
		$$ \texttt{tol}^n = \frac{| \|\partial_{\displacement}\vec R(\displacement^{n-1})[\delta\displacement^n] - \vec R(\displacement^{n-1})\| - \|\vec R(\displacement^{n-1})\| |} {\|\vec R(\displacement^{n-1})\|}. $$
		The choice of the norm $\|\cdot\|$ is arbitrary, so we consider in what follows the $\ell^2$ norm. This of course makes sense only in the discrete setting, meaning that for each residual vector $\vec r$ we consider the norm $\| \vec r\| = \left(\sum_i r_i^2 \right)^{1/2}$. The resulting scheme guarantees superlinear convergence \cite{dembo1982inexact}, which is of course worse than the quadratic convergence of a standard Newton-Krylov scheme, but gives an overall reduced complexity as it avoids oversolving the linearized problem.
	\item[BFGS method.] This is a Quasi-Newton method, and although it was initially devised as a minimization procedure \cite{wright1999numerical}, it can also be adapted to be used as a root-finding algorithms \cite{dennis1996numerical}. For this, we look at the minimization principle from which \eqref{eq:momentum strong} in quasi-static form is obtained:
		\begin{equation}\label{eq:minimization}
			\min_{\displacement} \Pi(\displacement)\coloneqq \int_\Omega \Psi(\ten F) - \vec f\cdot \displacement\,d\vec x .
		\end{equation}
A quadratic approximation of this problem yields
		$$\Pi(\displacement) \approx \Pi(\bar\displacement) + \partial_{\displacement}\Pi(\bar\displacement)[\displacement - \bar\displacement] + \frac 1 2 \partial^2_{\displacement} \Pi(\bar\displacement)[\displacement - \bar \displacement] ,$$
		and such scheme is indeed equivalent to \eqref{eq:newton} when applied as an iterative procedure with an exact Hessian. This method requires an initial approximation of the Hessian $\mat B^{-1,0} \approx [\partial^2_{\displacement} \Pi(\displacement^0)]^{-1}$, which is then enriched at each iteration with a rank two perturbation that includes curvature information given by
		$$ \mat B^{-1, k+1} = (\mat I - \rho_k \vec s^k \otimes \vec y^{k})\mat B^{-1, k}(\mat I - \rho_k \vec y^k \otimes \vec s^k) + \rho_k \vec s^k \otimes \vec s^{k}, $$
		where $\vec s^k=\vec x^{k+1} - \vec x^k$, $\vec y^k = \vec F(\vec x^{k+1}) - \vec F(\vec x^k)$, and $\rho^k= 1 / \langle \vec s^k, \vec y^k\rangle$ with $\langle\vec a,\vec b\rangle = \sum_ia_i b_i$. The approximate Hessian is then used to compute the next iteration as in Equation \eqref{eq:newton}:
			$$ \delta \displacement^k = -\mat B^{-1,k+1}\partial_\displacement \Pi(\displacement^{k-1}).$$
		The action of $\mat B^{-1,k+1}$ is implemented by means of a two-level recursion that allows for a limited memory implementation \cite{nocedal1980updating}. This method, as the inexact Newton-Krylov, yields superlinear convergence, and has the additional advantage of requiring only one assembly of the Hessian matrix (or any other initial matrix). The initial approximation of the Hessian is critical for the convergence and performance of this method, and we have indeed observed that simpler approaches do not yield satisfactory results for this method in the context of nonlinear mechanics. Motivated by this, we leverage the preconditioners obtained from the initial matrix $\partial_\displacement^2\Pi(\displacement^0)$ to obtain better approximations of the Hessian. If we denote by $\mat P$ the preconditioner arising from the initial Hessian, we consider the following approximations:
		\begin{itemize}
			\item $\mat B^{-1,0}=\text{action of $\mat P$ obtained from Hessian matrix} $\hfill (BFGS-preonly)
			\item $\mat B^{-1,0} = \text{inexact Krylov solver, large ($\approx 10^{-2}$) relative tolerance}$\hfill (inexact BFGS)
			\item $\mat B^{-1,0} = \text{quasi-exact Krylov solver, small ($\approx 10^{-6}$) relative tolerance}$\hfill(quasi-exact BFGS)
		\end{itemize}
		We have added on the right the names corresponding to Figure \ref{fig:exactness taxonomy}. We note that the quasi-exact and inexact versions are not covered by the theory, as the action of $\mat B^{-1,0}$ changes at each iteration. This happens because the number of linear iterations required at each nonlinear iteration varies according to the tolerance. We highlight that the use of a fixed tolerance in the inexact scenario is considered for simplicity, as an equivalent Einstat-Walker type of adaptive tolerance could be considered. This is an interesting alternative and by no means a trivial one to set up, so we leave this for future work.
\end{description}
We show a classification of these methods according to the exactness of the linear solver in Figure \ref{fig:exactness taxonomy}, where we have fixed three degrees of exactness: preconditioner only (preonly\footnote{The name is the same as the PETSc option to by-pass the iterative solver.}), inexact and quasi-exact. The two families present different types of performance: the Newton family relies on a low iteration count, where the level of inexactness relaxes the linear solver at each iteration, at the cost of some additional Newton steps. The quasi-Newton family instead yields a larger number of iterations, but each iteration is much cheaper than a Newton iteration, as it does not requires the Jacobian matrix to be reassembled. The level of exactness in this case reduces the iteration count, at the cost of making the iterations more expensive. The extreme cases Newton-preonly and quasi-exact BFGS have shown to be noncompetitive in our experiments, so we do not consider them in what follows. Throughout the remaining parts of this work, whenever no confusion arises, we may refer to Newton-Krylov and inexact Newton-Krylov methods as simply Newton methods. Throughout this manuscript, we will denote the Newton-Krylov, inexact Newton-Krylov, BFGS-preonly (or simply BFGS) and inexact-BFGS methods as \newton, \inexactnewton, \bfgs\, and \inexactbfgs\, respectively, as shown in Figure \ref{fig:exactness taxonomy}. We also show in Appendix \ref{appendix:petsc} the commands used to run each of the four methods under consideration. 

\begin{figure}
	\newcommand{\mhei}{0.4}
	\newcommand{\hei}{-0.8}
	\newcommand{\heii}{-1.6}
	\centering
		\begin{tikzpicture}
			\node at (0,0) (ex)  {Less inexact };
			\node at (0.6\textwidth,0) (inex){ More inexact};
			\node at (0.3\textwidth,\mhei) (tol) {smaller $\leftarrow$ Relative tolerance of linear solver $\rightarrow$ larger};
			\draw [stealth-stealth] (ex.east) -- (inex.west);
			\node at (0,\hei) (newton1) {Newton-Krylov (\newton)};
			\node at (0.3\textwidth,\hei) (newton2) {inexact Newton-Krylov (\inexactnewton)};
			\node at (0.6\textwidth,\hei) (newton3) {\xcancel{Newton-preonly}};
			\node at (0,\heii) (bfgs1) {\xcancel{Quasi-exact BFGS}};
			\node at (0.3\textwidth,\heii) (bfgs2) {inexact-BFGS (\inexactbfgs)};
			\node at (0.6\textwidth,\heii) (bfgs3) {BFGS-preonly (\bfgs)};
		\end{tikzpicture}
		\caption{Description of exactness level in the linear solver and nonlinear solvers used in the tests, where we denote "preconditioner only" as "preonly". Newton-preonly and quasi-exact BFGS have been discarded from our experiments because they are not competitive: the first requires too many matrix assembles, the second spends too much time at each iteration.}
		\label{fig:exactness taxonomy}
\end{figure}
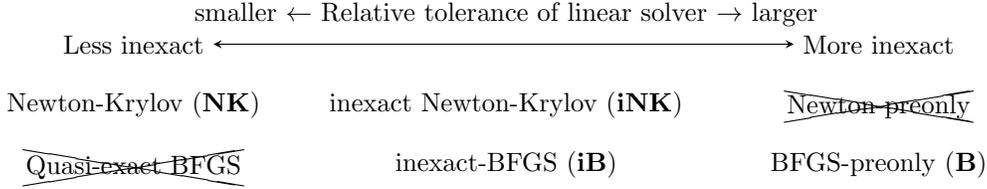

\section{Test problems and solvers details}\label{section:tests}
\begin{figure}
	\centering
	\begin{subfigure}{0.325\textwidth}
		\centering
		\includegraphics[height=5cm]{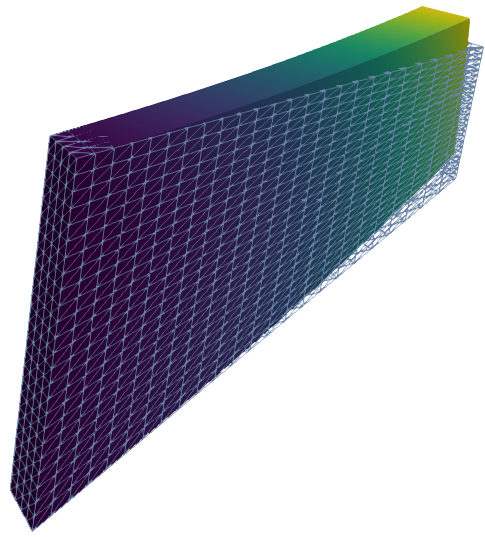}
		\caption{Cook test.}
	\end{subfigure}
	\begin{subfigure}{0.325\textwidth}
		\centering
		\includegraphics[height=5cm]{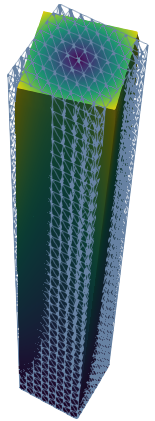}
		\caption{Twist test.}
	\end{subfigure}
	\begin{subfigure}{0.325\textwidth}
		\centering
		\includegraphics[height=5cm]{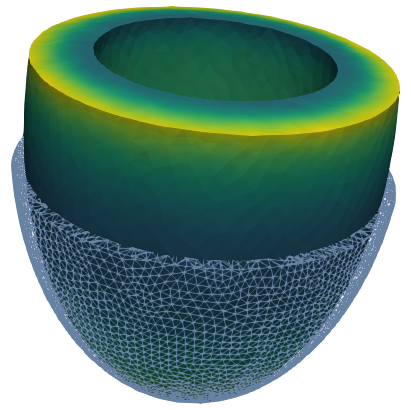}
		\caption{Heartbeat test.}
	\end{subfigure}
	\caption{Solutions of the benchmark tests used. Color is scaled according to displacement norm $|\displacement|$, and the reference configuration is depicted with its mesh. In (c), the mesh is cut roughly at half the height of the geometry.}
	\label{fig:solutions}
\end{figure}

In this section, we describe the test cases we use to compare the performance of each method. The scope of the heartbeat test is that of providing a more realistic scenario in which the insight acquired in the other two tests will be transfered. Details on the implementation will be given on each problem.
\begin{description}
	\item[Cook test.] This is a static almost incompressible mechanics benchmark \cite{Pelteret_McBride_2016}. We use homogeneous Dirichlet conditions on $\{x=0\}$, a vertical traction given by $\vec t = (0, \tau, 0)$ with $\tau=10^6$ in deformed configuration on $\{x=L\}$, a homogeneous distributed load $\vec f=\vec 0$, and homogeneous Neumann conditions on the remaining parts of the boundary. The constitutive modeling is given by a almost incompressible Neo-Hookean material:
		$$\Psi_\text{Cook}(\ten F) = C_1\left(\tr{\left(\bar{\ten F}^T \bar{\ten F}\right)} - 3\right) + k\left([\det{(\ten F)}]^2 - 1 - 2\log \left(\det \ten F\right)\right), $$
		where $C_1=\frac 1 2 \mu$, $\mu = 8.194\cdot10^7$, $k=\lambda + \frac 2 3 \mu$, $\lambda = 2\mu\nu/(1 - 2\nu)$, and $\nu =0.3$. The Krylov solver was configured with a GMRES linear solver without restart and with a modified Gram-Schmidt procedure which is more robust \cite{golub1996matrix}, preconditioned with HYPRE-BoomerAMG \cite{falgout2002hypre}, using its default configuration and setting the matrix block size to 3 in PETSc to improve the efficiency of the preconditioner. The absolute tolerance was set to $10^{-8}$ for the exact solver, and to $0$ for the inexact ones, whereas the relative tolerance was set to $10^{-6}$ for the exact solver and to $10^{-2}$ for the inexact ones. We show the solution in Figure \ref{fig:solutions} (a).

	\item[Twist test.] This is a static incompressible mechanics problem \cite{bonet2015computational}, where we use the inf-sup stable stable finite elements $\mathbb P_2-\mathbb P_0$ for the approximation of the displacement and pressure. Boundary conditions are given by a homogeneous on $\{z=0\}$, a $\pi/6$ rotation in the $x$ and $y$ components on $\{z=L_z\}$, homogeneous data $\vec f=\vec t=\vec 0$, and homogeneous Neumann conditions on the remaining parts of the boundary. The constitutive modelling is given by the following polyconvex potential:
		$$ \Psi_\text{Twist}(\ten F,p) = \alpha_p(\ten F:\ten F - 3) + \beta_p (\cof \ten F:\cof \ten F - 3) - (4\beta_s + 2\alpha_s)\log\left(\ten F\right) - p(\det(\ten F) - 1),$$
		where $\alpha_p = \beta_p = \alpha_s = \beta_s = 9000$ and $\cof \ten F\coloneqq \det(\ten F) \ten F^{-T}$. The logarithmic term with $\det(\ten F)$ is required because of the weak imposition of the incompressibility constraint \cite{auricchio2013approximation}. The preconditioning of a saddle point problem is more challenging, so to obtain a more robust performance we used a lower Schur complement preconditioner based on a $(\displacement, p)$ field split \cite{brown2012composable}, we provide details on the Schur complement preconditioner in Appendix \ref{appendix:schur}. The displacement block uses an AMG preconditioner, the Schur complement block instead uses the SIMPLE preconditioner from fluid mechanics \cite{elman2006block}, meaning that the inverse of the displacement block is approximated by the inverse of its diagonal, and for the resulting block we used a block Jacobi preconditioner, which we observed to suffice. Two additional comments are in place: (i) to calibrate the effectiveness of each solver in the Schur complement, we started from using a direct solver in each block to guarantee convergence in at most two iterations \cite{mandel1990block}, and started relaxing the blocks from there. The difficulty of the problem is predominantly the displacement block, which can be appreciated by the simple preconditioner used for the Schur complement block; (ii) the pressure block presents much less degrees of freedom (DoFs), so we used a telescopic approach that creates an MPI subcommunicator with a fixed $25\%$ of the original processes to avoid excessive communication. We show the solution in Figure \ref{fig:solutions} (b).

	\item[Heartbeat test.] This problem represents a real-case study and is indeed our case of interest. It models the contraction of a human left ventricle in an idealized geometry given by a prolate ellipsoid, with a pointwise set of coordinates representing the muscle fiber orientation $(\vec f_0, \vec s_0, \vec t_0)$ obtained through rule based methods \cite{Bayer20122243}. We consider a Guccione hyperelastic potential \cite{Guccione1991}, given by
	\begin{align*}
		 \Psi(\ten F) &= \frac C 2 \left(e^{Q(\ten F)} - 1\right) + \frac B 2 (J-1)\log J, \\
		 Q(\ten F) &= b_{ff}E_{ff}^2 + b_{ss}E_{ss}^2 + b_{nn}E_{nn}^2 + b_{fs}(E_{fs}^2 + E_{sf}^2) + b_{fn}(E_{fn}^2 + E_{nf}^2) + b_{sn}(E_{sn}^2 + E_{ns}^2),
		\end{align*}
		where $\ten E = \frac 1 2 \left( \ten F^T\ten F - \ten I\right)$ and $E_{ab} = \ten E\vec a\cdot \vec b$ for $a,b\in\{f,s,n\}$, which are the components of $\ten E$ in the fiber-induced frame of reference $\{\vec f, \vec s, \vec n\}$. See \cite{Guccione1991} for reference values of the related parameters. The Piola stress tensor is enriched with a fiber-wise force known as active stress that models the contraction of the muscle cells (cardiomyocites):
		$$ \ten P(\ten F) = \frac{\partial \Psi}{\partial \ten F}(\ten F) + \ten P_\text{act}(\ten F), \quad \ten P_\text{act}(\ten F) = \gamma(t)\frac{(\ten F \vec f)\otimes \vec f}{|\ten F \vec f|}, $$
		where $\gamma$ is a given function known as the activation function. For simplicity we consider an analytical activation given by 
		$$\gamma(t) = C_\text{PA}\max\{\sin (2\pi t/T), 0\} ,$$
		where the peak activation constant is $C_\text{PA}=10^4$ and the period is $T=0.8$. More details on the generation of the fibers and the activation function can be found in \cite{Bayer20122243} and \cite{regazzoni2020cardiac} respectively.  Boundary conditions are given by homogeneous Neumann on the endocardium $\partial\Omega_\text{endo}$ and the base $\partial\Omega_\text{base}$, whereas the epicardium $\partial\Omega_\text{epi}$ considers a Robin condition that models the friction of the epicardium with the pericardium \cite{Usyk2002} and is given by
		\begin{equation}\label{eq:bc epi}
			\ten P(\ten F)\vec n = \coloneqq - (\vec n\otimes \vec n)\left(K_\perp^\text{epi}\displacement + C_\perp^\text{epi}\dot\displacement\right) - (\ten I-\vec n \otimes \vec n)\left(K_\|^\text{epi}\displacement + C_\|^\text{epi}\dot\displacement\right)=\vec 0\quad\text{ on }\partial\Omega_\text{epi}.
		\end{equation}
		More information of the physical meaning of these boundaries and the motivation for these choices can be found in \cite{quarteroni2017integrated}. The Krylov solver in this case is identical to the one used in the Cook test. We show the solution at $t=0.2$ in Figure \ref{fig:solutions} (c).
\end{description}

\section{Theoretical results and solvers expected behavior}\label{section:theory}
The purpose of this section is to show that the hypotheses from the convergence analysis, despite not being directly applicable to this problem, provide fundamental insight to guide the choice of the nonlinear solver. To state the results we consider an arbitrary scalar function $f:\mathbb R^m \to \mathbb R$ to be minimized and an initial point $\vec x^0$. The results and their proofs can be found in \cite{wright1999numerical}. First we state the convergence result of Newton's method, where we remark that the inexact variant requires no additional hypotheses.

\begin{theorem}[Newton methods]
	 Assume that there exists a point $\vec x^*$ such that $\nabla f(\vec x^*)=0$, $\nabla^2 f$ is Lipschitz continuous in a neighborhood of $\vec x^*$ and $\nabla^2 f(\vec x^*)$ is positive definite. Then, $\vec x^*$ is a local minimizer and, if the initial guess $\vec x^0$ is sufficiently close to $\vec x^*$, the iterates obtained by Newton's method converge quadratically to $\vec x^*$.
\end{theorem}
Now we state the general convergence theorem for the BFGS method.
\begin{theorem}[Quasi-Newton methods]
	Assume that $f$ is twice continuously differentiable and that $\nabla^2 f$ is Lipschitz continuous at the minimizer $\vec x^*$. Assume also that the level
		$$ \mathcal L = \{\vec x \in \mathbb R^m: f(\vec x)\leq f(\vec x^0)\}, $$
	is convex, that $\nabla^2 f$ is positive definite on $\mathcal L$ and that the initial matrix $\mat B^0$ is positive definite. Then, the sequence of iterates obtained by the BFGS method converges superlinearly to $\vec x^*$.
\end{theorem}
We make the following observations, which are fundamental to understand the performance of these methods.
\begin{itemize}
	\item Both methods  require the initial point to be sufficiently close. This can be seen in the BFGS method through the convexity of $\mathcal L$.
	\item The BFGS method is better suited for convex problems, given that the Hessian is required to be definite positive not only at the minimizer as Newton, but also in the entire level set $\mathcal L$.
	\item Inexact Newton-Krylov convergence holds under the same hypotheses of Newton's method, plus some hypotheses on the relative tolerances $\eta_k$ that guarantee the superlinear convergence. We consider tolerances that satisfy such hypotheses \cite{eisenstat1994globally}, and note that they can sometimes yield global convergence, unlike Newton's method. This makes inexact Newton-Krylov methods attractive as they can be potentially more robust than a classic Newton method.
	\item The initial matrix approximation $\mat B^0$ is variable in the inexact BFGS method, meaning that its performance might not necessarily be better, nor more robust, than the case in which only the action of the preconditioner is considered. This is due to this case not being covered by the theory, so that in principle there could be additional hypotheses on the action of a variable $\mat B^0$.
	\item Notes for different types of problems:
		\begin{itemize}
			\item {\bf Cook test:} In this case, the problem being solved is polyconvex, but if a sufficiently close initial guess is considered, we can possibly fall in an attraction basin. This means that both Newton and BFGS methods should convergence, but both BFGS might possibly be less robust due to the convexity requirement of $\mathcal L$.
			\item {\bf Twist test:} This is a saddle point problem, so BFGS is not guaranteed to converge even if a sufficiently good initial approximation is considered. This can be seen from the positive definiteness hypothesis of the Hessian on $\mathcal L$.
			\item {\bf Heartbeat test:} The inertia term in elastodynamics results in a convex contribution to the variational principle associated to the problem, so we expect this case to be easier to solve in spite of the more complex nonlinearities involved. Instead, the nonlinearities should impact the number of nonlinear iterations. 
		\end{itemize}
\end{itemize}

\section{Numerical results}\label{section:results}
In this section, we compare the performance of the solvers introduced in Sec. \ref{section:solvers}, for all tests under consideration. The scope of these tests is to assess the sensitivity of the methods with respect to the problem size, some of the problem parameters, and the number of computational cores. For the heartbeat problem, we consider four different meshes with the degrees of freedom shown in Table \ref{table:heart meshes dofs}. We note that meshes 3 and 4 are refinements of meshes 1 and 2 respectively. Note as well that since BFGS-preonly does not solve any linear system, we do not report the linear iterations for this method. All models are implemented using the FEniCS library \cite{alnaes2015fenics}, and all interfaces with the underlying PETSc library are done by using the \texttt{petsc4py} interface \cite{petsc-user-ref}. The computations were performed on the EOS, INDACO and Galileo100 supercomputers.

\begin{table}[h!]
	\centering
	\begin{tabular}{c | r  r}
	\toprule Mesh name & DoFs \Pone & DoFs \Ptwo \\
	\midrule Mesh 1 & 9375 & 60081 \\
	Mesh 2 & 20709 & 149511 \\
	Mesh 3 & 60081 & 421191 \\
	Mesh 4 & 149511 & 1123515 \\ \bottomrule
	\end{tabular}
	\caption{Heartbeat test: Degrees of freedom yielded by each of the meshes under consideration.}
	\label{table:heart meshes dofs}
\end{table}

\begin{table}
\begin{center}
\begin{tabular}{r|rrr|rrr|rrr|rr}
\toprule
\multicolumn{12}{c}{\Pone \ discretization} \\
\toprule
DoFs            & \multicolumn{3}{c|}{\newton}        & \multicolumn{3}{c|}{\inexactnewton}       & \multicolumn{3}{c|}{\inexactbfgs}       & \multicolumn{2}{c}{\bfgs}\\
                & nit   & lit   & $T_{sol}$       & nit   & lit   & $T_{sol}$       & nit   & lit   & $T_{sol}$       & nit   & $T_{sol}$\\
\midrule
9375            & 4     & 11.8  & 0.7           & 5     & 9.4   & 0.8           & 6     & 8.7   & 0.4           & 112   & 4.0 \\
64827           & 5     & 11.6  & 5.4           & 5     & 11.6  & 5.5           & 6     & 9.7   & 2.7           & 178   & 46.5 \\
207831          & 6     & 10.3  & 23.2          & 5     & 10.6  & 20.3          & 7     & 10.0  & 12.7          & 59    & 17.2 \\
479859          & 5     & 13.0  & 51.0          & 5     & 12.6  & 52.0          & 7     & 9.6   & 30.5          & 46    & 31.8 \\
922383          & 5     & 13.4  & 105.0         & 5     & 13.0  & 103.2         & 7     & 11.1  & 69.5          & 43    & 52.7 \\
1576875         & 5     & 13.8  & 178.6         & 5     & 13.4  & 183.2         & 7     & 10.1  & 115.7         & 48    & 112.7 \\
\bottomrule
\toprule
\multicolumn{12}{c}{\Ptwo \ discretization} \\
\toprule
DoFs            & \multicolumn{3}{c|}{\newton}        & \multicolumn{3}{c|}{\inexactnewton}       & \multicolumn{3}{c|}{\inexactbfgs}        & \multicolumn{2}{c}{\bfgs}\\
                & nit   & lit   & $T_{sol}$     & nit   & lit   & $T_{sol}$     & nit   & lit   & $T_{sol}$     & nit   & $T_{sol}$\\
\midrule
9375            & 5     & 15.8  & 1.6           & 5     & 18.6  & 1.6           & 7     & 14.6  & 0.8           & 89    & 1.7 \\
64827           & 4     & 20.2  & 11.5          & 5     & 17.4  & 13.6          & 7     & 14.1  & 7.8           & 203   & 65.3 \\
207831          & 5     & 16.8  & 50.5          & 5     & 17.2  & 50.7          & 7     & 14.0  & 31.8          & 65    & 31.1 \\
479859          & 4     & 20.0  & 102.2         & 5     & 18.2  & 118.0         & 7     & 14.1  & 79.0          & 326   & 1421.0\\
922383          & 4     & 20.8  & 198.7         & 5     & 17.6  & 226.0         & 8     & 13.3  & 168.8         & 332   & 2692.2 \\
1576875         & 4     & 20.5  & 358.3         & 5     & 17.0  & 406.6         & 7     & 13.3  & 268.8         & 64    & 240.3 \\
\bottomrule
\end{tabular}
\caption{Sensitivity with respect to problem size, Cook test. nit$\,\coloneqq$ nonlinear iterations, lit$\,\coloneqq$ average linear iterations per nonlinear iteration, $T_{sol}\coloneqq$ CPU time in seconds.}
\label{tab_perf_cook}
\end{center}
\end{table}

\begin{table}
\begin{center}
\begin{tabular}{r|rrr|rrr}
\toprule
DoFs            & \multicolumn{3}{c}{\newton}         & \multicolumn{3}{c}{\inexactnewton} \\
                & nit   & lit   & $T_{sol}$       & nit   & lit   & $T_{sol}$ \\
\midrule
131163          & 5     & 551.8 & 836.7         & 15    & 162.8 & 980.2 \\
223347          & 5     & 569.2 & 2045.5        & 15    & 170.1 & 1989.2 \\
350955          & 6     & 699.7 & 5005.9        & 15    & 215.1 & 4100.7 \\
519747          & 6     & 834.7 & 9332.2        & 16    & 233.7 & 7347.4 \\
\bottomrule
\end{tabular}
\caption{Sensitivity with respect to problem size, Twist test. nit$\,\coloneqq$ nonlinear iterations, lit$\,\coloneqq$ average linear iterations per nonlinear iteration, $T_{sol}\coloneqq$ CPU time in seconds.}
\label{tab_perf_twist}
\end{center}
\end{table}

\begin{table}
\begin{center}
\begin{tabular}{r|rrr|rrr|rrr|rr}
\toprule
\multicolumn{12}{c}{\Pone \ discretization} \\
\toprule
DoFs            & \multicolumn{3}{c|}{\newton}        & \multicolumn{3}{c|}{\inexactnewton}       & \multicolumn{3}{c|}{\inexactbfgs}       & \multicolumn{2}{c}{\bfgs}\\
                & nit   & lit   & $T_{sol}$       & nit   & lit   & $T_{sol}$       & nit   & lit   & $T_{sol}$       & nit   & $T_{sol}$\\
\midrule
9375            & 3.6   & 7.6   & 0.4           & 4.8   & 5.6   & 0.5           & 7.9   & 3.0   & 0.2           & 24.6  & 0.3 \\
20709           & 3.6   & 7.1   & 1.1           & 4.6   & 5.0   & 1.3           & 8.1   & 2.8   & 0.6           & 22.1  & 0.8 \\
60081           & 3.6   & 8.3   & 3.0           & 5.0   & 5.8   & 3.7           & 8.3   & 3.0   & 1.8           & 25.9  & 2.3 \\
149511          & 3.7   & 8.0   & 9.9           & 4.9   & 5.7   & 11.9          & 8.1   & 3.2   & 6.2           & 24.4  & 7.2 \\
\bottomrule
\toprule
\multicolumn{12}{c}{\Ptwo \ discretization} \\
\toprule
DoFs            & \multicolumn{3}{c|}{\newton}        & \multicolumn{3}{c|}{\inexactnewton}       & \multicolumn{3}{c|}{\inexactbfgs}        & \multicolumn{2}{c}{\bfgs}\\
                & nit   & lit   & $T_{sol}$     & nit   & lit   & $T_{sol}$     & nit   & lit   & $T_{sol}$     & nit   & $T_{sol}$\\
\midrule
60081           & 3.7   & 10.4  & 4.4           & 5.0   & 7.0   & 5.1           & 9.6   & 3.9   & 2.6           & 32.7  & 2.5 \\
149511          & 3.7   & 8.9   & 13.2          & 5.0   & 6.0   & 15.4          & 9.4   & 3.3   & 8.2           & 26.5  & 6.9 \\
421191          & 3.8   & 10.8  & 41.4          & 5.2   & 7.1   & 45.7          &       & F     &               & 33.0  & 23.5 \\
1123515         & 3.8   & 11.0  & 125.5         & 5.1   & 7.2   & 137.0         &       & F     &               & 32.2  & 70.9 \\
\bottomrule
\end{tabular}
\caption{Sensitivity with respect to problem size, Heartbeat test. nit$\,\coloneqq$ nonlinear iterations, lit$\,\coloneqq$ average linear iterations per nonlinear iteration, $T_{sol}\coloneqq$ CPU time in seconds.}
\label{tab_perf_heart}
\end{center}
\end{table}

\begin{figure}[ht]
\begin{center}
\input{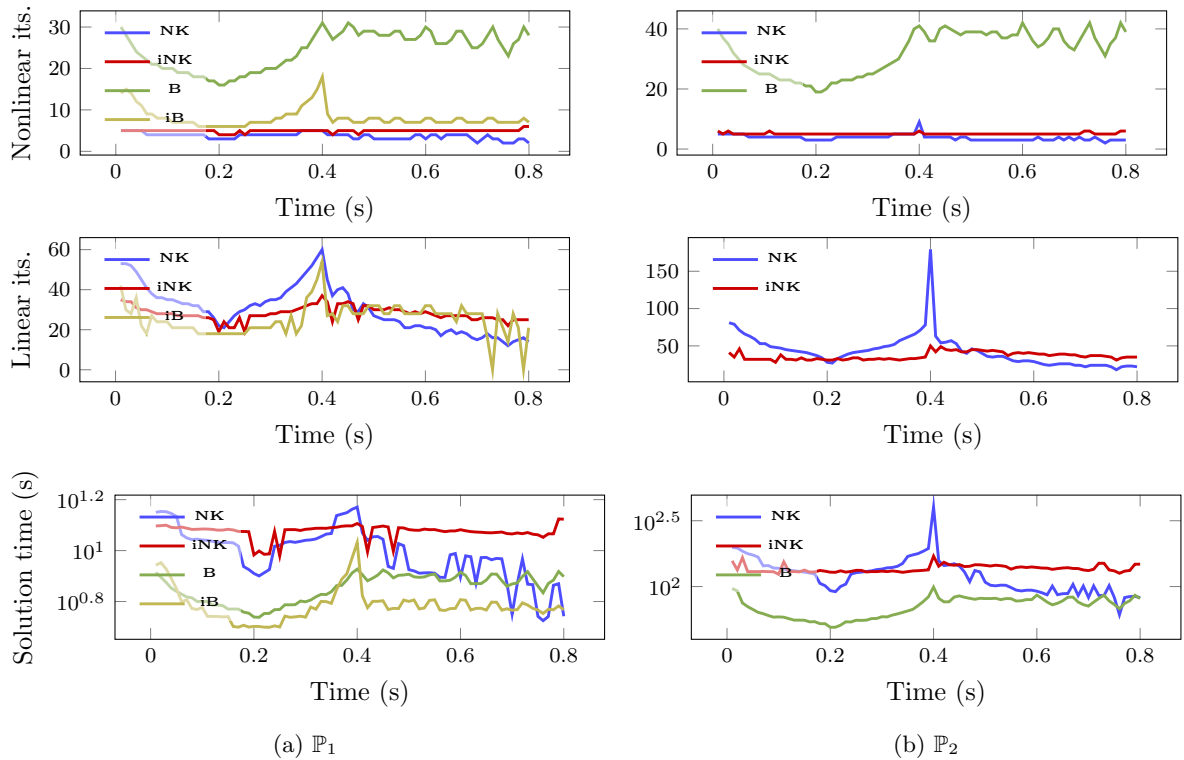}
\caption{Sensitivity with respect to problem size, Heartbeat test. Time evolution of nonlinear iterations (nit), average linear iterations per nonlinear iteration (lit) and CPU time in seconds for \Pone \ and \Ptwo \ discretizations in mesh 4.}
\label{fig_perf_heart}
\end{center}
\end{figure}

\subsection{Sensitivity with respect to problem size}
To compare the performance of all methods, we report the total nonlinear iterations counts, the average Krylov iterations per nonlinear step and the solution time as we increase the number of degrees of freedom.
\begin{description}
        \item[Cook test.] We present the nonlinear iterations, linear iterations and CPU time in Table \ref{tab_perf_cook}, varying the dimension of the problem from about 10 thousand to 1.5 million degrees of freedom. We first observe that notably all methods exhibit a robust behavior with respect to the degrees of freedom, except for the BFGS-preonly with second order finite elements. Indeed, the total and average linear iterations reported are roughly constant, but in general the inexact-BFGS yields lower average linear iterations. All methods exhibit an approximately linear increase of CPU time with respect to the degrees of freedom, with the exception of BFGS-preonly with second order elements. In general, BFGS methods outperform Newton methods in terms of the execution time, but attention must be given to the second order scenario. The ratio between the CPU times of the best method (BFGS-preonly) and Newton's method in the case with the largest amount of DoFs is 0.63 for \Pone and 0.67 for \Ptwo elements.

        \item[Twist test.] We report the nonlinear iterations, average linear iterations per nonlinear step, and CPU times in Table \ref{tab_perf_twist}. We observe that BFGS methods do not converge for this problem, so we can not report results for them. This holds not only for both methods presented (BFGS, inexact-BFGS), as we also tried the BFGS-exact variant described in Section \ref{section:solvers}, which was ineffective as well.
        Both Newton-Krylov and inexact Newton-Krylov methods present a very mild increase in the nonlinear iterations, but the inexact variant presents a plateau of average linear iterations, whereas the Newton-Krylov method presents a monotonic increase in the average linear iterations. This results in the inexact method being a faster method, and the difference in CPU times increases with the problem size. The ratio between the CPU times of the inexact Newton-Krylov method and Newton-Krylov method in the case with the largest amount of DoFs is 0.79.

        \item[Heartbeat test.] We report the nonlinear iterations, average linear iterations per nonlinear iteration, and CPU time for first and second order finite elements averaged during the first timestep in Table \ref{tab_perf_heart}  and their time evolution in Figure \ref{fig_perf_heart}. We note that all methods are robust with respect to the degrees of freedom, except the inexact-BFGS, as it does not converge for fine meshes when using second order elements. Both Newton methods present similar nonlinear iteration counts, and all methods have a similar average linear iteration count. This results in the inexact-BFGS having a lower average linear iteration count overall as it requires more nonlinear iterations. This has an impact in the overall solution time, where BFGS methods perform much better than Newton methods. The ratio between the CPU times of the best method (inexact-BFGS in \Pone, BFGS-preonly in \Ptwo) and the Newton-Krylov method is 0.63 for \Pone and 0.57 for \Ptwo elements.
\end{description}

\begin{figure}
  \begin{subfigure}{\textwidth}
  \begin{subfigure}[b]{0.325\textwidth}
      \begin{semilogxplot}{$\tau$}{Nonlinear its.}
      \addplot+[fill opacity=0.1, color=\newtoncolor] table [x=rob-param, y=nl-its, col sep=comma] {data/result-rob-cook-Newton-P1.csv};
      \addplot+[fill opacity=0.1, color=\inexactnewtoncolor] table [x=rob-param, y=nl-its, col sep=comma] {data/result-rob-cook-Newton-inexact-P1.csv};
      \addplot+[fill opacity=0.1, color=\bfgscolor] table [x=rob-param, y=nl-its, col sep=comma] {data/result-rob-cook-BFGS-P1.csv};
      \addplot+[fill opacity=0.1, color=\inexactbfgscolor] table [x=rob-param, y=nl-its, col sep=comma] {data/result-rob-cook-BFGS-inexact-P1.csv};
      \legend {\newton, \inexactnewton, \bfgs, \inexactbfgs}
    \end{semilogxplot}
  \end{subfigure}
  \begin{subfigure}[b]{0.325\textwidth}
      \begin{plot}{$\tau$}{Total linear its.}
      \addplot+[fill opacity=0.1, color=\newtoncolor] table [x=rob-param, y=tot-krylov, col sep=comma] {data/result-rob-cook-Newton-P1.csv};
      \addplot+[fill opacity=0.1, color=\inexactnewtoncolor] table [x=rob-param, y=tot-krylov, col sep=comma] {data/result-rob-cook-Newton-inexact-P1.csv};
      \addplot+[fill opacity=0.1, color=\inexactbfgscolor] table [x=rob-param, y=tot-krylov, col sep=comma] {data/result-rob-cook-BFGS-inexact-P1.csv};
      \legend {\newton, \inexactnewton, \inexactbfgs}
    \end{plot}
  \end{subfigure}
  \begin{subfigure}[b]{0.325\textwidth}
      \begin{plot}{$\tau$}{Average linear its.}
      \addplot+[fill opacity=0.1, color=\newtoncolor] table [x=rob-param, y=avg-krylov, col sep=comma] {data/result-rob-cook-Newton-P1.csv};
      \addplot+[fill opacity=0.1, color=\inexactnewtoncolor] table [x=rob-param, y=avg-krylov, col sep=comma] {data/result-rob-cook-Newton-inexact-P1.csv};
      \addplot+[fill opacity=0.1, color=\inexactbfgscolor] table [x=rob-param, y=avg-krylov, col sep=comma] {data/result-rob-cook-BFGS-inexact-P1.csv};
      \legend {\newton, \inexactnewton, \inexactbfgs}
    \end{plot}
  \end{subfigure}
  \caption{\Pone}
  \end{subfigure}

  \begin{subfigure}{\textwidth}
  \begin{subfigure}[b]{0.325\textwidth}
      \begin{semilogxplot}{$\tau$}{Nonlinear its.}
      \addplot+[fill opacity=0.1, color=\newtoncolor] table [x=rob-param, y=nl-its, col sep=comma] {data/result-rob-cook-Newton-P2.csv};
      \addplot+[fill opacity=0.1, color=\inexactnewtoncolor] table [x=rob-param, y=nl-its, col sep=comma] {data/result-rob-cook-Newton-inexact-P2.csv};
      \addplot+[fill opacity=0.1, color=\bfgscolor] table [x=rob-param, y=nl-its, col sep=comma] {data/result-rob-cook-BFGS-P2.csv};
      \addplot+[fill opacity=0.1, color=\inexactbfgscolor] table [x=rob-param, y=nl-its, col sep=comma] {data/result-rob-cook-BFGS-inexact-P2.csv};
      \legend {\newton, \inexactnewton, \bfgs, \inexactbfgs}
    \end{semilogxplot}
  \end{subfigure}
  \begin{subfigure}[b]{0.325\textwidth}
      \begin{plot}{$\tau$}{Total linear its.}
      \addplot+[fill opacity=0.1, color=\newtoncolor] table [x=rob-param, y=tot-krylov, col sep=comma] {data/result-rob-cook-Newton-P2.csv};
      \addplot+[fill opacity=0.1, color=\inexactnewtoncolor] table [x=rob-param, y=tot-krylov, col sep=comma] {data/result-rob-cook-Newton-inexact-P2.csv};
      \addplot+[fill opacity=0.1, color=\inexactbfgscolor] table [x=rob-param, y=tot-krylov, col sep=comma] {data/result-rob-cook-BFGS-inexact-P2.csv};
      \legend {\newton, \inexactnewton, \inexactbfgs}
    \end{plot}
  \end{subfigure}
  \begin{subfigure}[b]{0.325\textwidth}
      \begin{plot}{$\tau$}{Average linear its.}
      \addplot+[fill opacity=0.1, color=\newtoncolor] table [x=rob-param, y=avg-krylov, col sep=comma] {data/result-rob-cook-Newton-P2.csv};
      \addplot+[fill opacity=0.1, color=\inexactnewtoncolor] table [x=rob-param, y=avg-krylov, col sep=comma] {data/result-rob-cook-Newton-inexact-P2.csv};
      \addplot+[fill opacity=0.1, color=\inexactbfgscolor] table [x=rob-param, y=avg-krylov, col sep=comma] {data/result-rob-cook-BFGS-inexact-P2.csv};
      \legend {\newton, \inexactnewton, \inexactbfgs}
    \end{plot}
  \end{subfigure}
  \caption{\Ptwo}
  \end{subfigure}
  \caption{Robustness with respect to the data, Cook test. Comparison of the nonlinear, total linear  and average linear iterations.}
  \label{fig:cook rob all}
\end{figure}

\begin{figure}
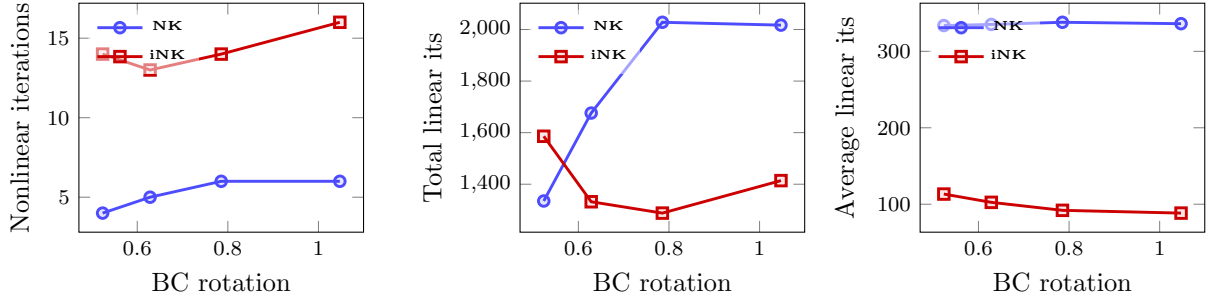

  \centering
  \begin{subfigure}[b]{0.325\textwidth}
      \begin{plot}{BC rotation}{Nonlinear iterations}
      \addplot+[fill opacity=0.1, color=\newtoncolor] table [x=rob-param, y=nl-its, col sep=comma] {data/result-rob-twist-Newton.csv};
      \addplot+[fill opacity=0.1, color=\inexactnewtoncolor] table [x=rob-param, y=nl-its, col sep=comma] {data/result-rob-twist-Newton-inexact.csv};
      \legend {\newton, \inexactnewton}
    \end{plot}
  \end{subfigure}
  \begin{subfigure}[b]{0.325\textwidth}
      \begin{plot}{BC rotation}{Total linear its}
      \addplot+[fill opacity=0.1, color=\newtoncolor] table [x=rob-param, y=tot-krylov, col sep=comma] {data/result-rob-twist-Newton.csv};
      \addplot+[fill opacity=0.1, color=\inexactnewtoncolor] table [x=rob-param, y=tot-krylov, col sep=comma] {data/result-rob-twist-Newton-inexact.csv};
      \legend {\newton, \inexactnewton}
    \end{plot}
  \end{subfigure}
  \begin{subfigure}[b]{0.325\textwidth}
      \begin{plot}{BC rotation}{Average linear its}
      \addplot+[fill opacity=0.1, color=\newtoncolor] table [x=rob-param, y=avg-krylov, col sep=comma] {data/result-rob-twist-Newton.csv};
      \addplot+[fill opacity=0.1, color=\inexactnewtoncolor] table [x=rob-param, y=avg-krylov, col sep=comma] {data/result-rob-twist-Newton-inexact.csv};
      \legend {\newton, \inexactnewton}
    \end{plot}
  \end{subfigure}
  \caption{Robustness with respect to the data, Twist test. Nonlinear, total linear and average linear iterations with respect to the rotation angle.}
  \label{fig:twist robust all}
\end{figure}

\begin{figure}
  \begin{subfigure}{\textwidth}
  \begin{subfigure}[b]{0.325\textwidth}
      \begin{plot}{$C_{PA}$}{Nonlinear its.}
      \addplot+[fill opacity=0.1, color=\newtoncolor] table [x=rob-param, y=nl-its, col sep=comma] {data/result-rob-heart-Newton-P1.csv};
      \addplot+[fill opacity=0.1, color=\inexactnewtoncolor] table [x=rob-param, y=nl-its, col sep=comma] {data/result-rob-heart-Newton-inexact-P1.csv};
      \addplot+[fill opacity=0.1, color=\bfgscolor] table [x=rob-param, y=nl-its, col sep=comma] {data/result-rob-heart-BFGS-P1.csv};
      \addplot+[fill opacity=0.1, color=\inexactbfgscolor] table [x=rob-param, y=nl-its, col sep=comma] {data/result-rob-heart-BFGS-inexact-P1.csv};
      \legend {\newton, \inexactnewton, \bfgs, \inexactbfgs}
    \end{plot}
  \end{subfigure}
  \begin{subfigure}[b]{0.325\textwidth}
      \begin{plot}{$C_{PA}$}{Total linear its.}
      \addplot+[fill opacity=0.1, color=\newtoncolor] table [x=rob-param, y=tot-krylov, col sep=comma] {data/result-rob-heart-Newton-P1.csv};
      \addplot+[fill opacity=0.1, color=\inexactnewtoncolor] table [x=rob-param, y=tot-krylov, col sep=comma] {data/result-rob-heart-Newton-inexact-P1.csv};
      \addplot+[fill opacity=0.1, color=\inexactbfgscolor] table [x=rob-param, y=tot-krylov, col sep=comma] {data/result-rob-heart-BFGS-inexact-P1.csv};
      \legend {\newton, \inexactnewton, \inexactbfgs}
    \end{plot}
  \end{subfigure}
  \begin{subfigure}[b]{0.325\textwidth}
      \begin{plot}{$C_{PA}$}{Average linear its.}
      \addplot+[fill opacity=0.1, color=\newtoncolor] table [x=rob-param, y=avg-krylov, col sep=comma] {data/result-rob-heart-Newton-P1.csv};
      \addplot+[fill opacity=0.1, color=\inexactnewtoncolor] table [x=rob-param, y=avg-krylov, col sep=comma] {data/result-rob-heart-Newton-inexact-P1.csv};
      \addplot+[fill opacity=0.1, color=\inexactbfgscolor] table [x=rob-param, y=avg-krylov, col sep=comma] {data/result-rob-heart-BFGS-inexact-P1.csv};
      \legend {\newton, \inexactnewton, \inexactbfgs}
    \end{plot}
  \end{subfigure}
  \caption{\Pone}
  \end{subfigure}

  \begin{subfigure}{\textwidth}
  \begin{subfigure}[b]{0.325\textwidth}
      \begin{plot}{$C_{PA}$}{Nonlinear its.}
      \addplot+[fill opacity=0.1, color=\newtoncolor] table [x=rob-param, y=nl-its, col sep=comma] {data/result-rob-heart-Newton-P2.csv};
      \addplot+[fill opacity=0.1, color=\inexactnewtoncolor] table [x=rob-param, y=nl-its, col sep=comma] {data/result-rob-heart-Newton-inexact-P2.csv};
      \addplot+[fill opacity=0.1, color=\bfgscolor] table [x=rob-param, y=nl-its, col sep=comma] {data/result-rob-heart-BFGS-P2.csv};
      \addplot+[fill opacity=0.1, color=\inexactbfgscolor] table [x=rob-param, y=nl-its, col sep=comma] {data/result-rob-heart-BFGS-inexact-P2.csv};
      \legend {\newton, \inexactnewton, \bfgs, \inexactbfgs}
    \end{plot}
  \end{subfigure}
  \begin{subfigure}[b]{0.325\textwidth}
      \begin{plot}{$C_{PA}$}{Total linear its.}
      \addplot+[fill opacity=0.1, color=\newtoncolor] table [x=rob-param, y=tot-krylov, col sep=comma] {data/result-rob-heart-Newton-P2.csv};
      \addplot+[fill opacity=0.1, color=\inexactnewtoncolor] table [x=rob-param, y=tot-krylov, col sep=comma] {data/result-rob-heart-Newton-inexact-P2.csv};
      \addplot+[fill opacity=0.1, color=\inexactbfgscolor] table [x=rob-param, y=tot-krylov, col sep=comma] {data/result-rob-heart-BFGS-inexact-P2.csv};
      \legend {\newton, \inexactnewton, \inexactbfgs}
    \end{plot}
  \end{subfigure}
  \begin{subfigure}[b]{0.325\textwidth}
      \begin{plot}{$C_{PA}$}{Average linear its.}
      \addplot+[fill opacity=0.1, color=\newtoncolor] table [x=rob-param, y=avg-krylov, col sep=comma] {data/result-rob-heart-Newton-P2.csv};
      \addplot+[fill opacity=0.1, color=\inexactnewtoncolor] table [x=rob-param, y=avg-krylov, col sep=comma] {data/result-rob-heart-Newton-inexact-P2.csv};
      \addplot+[fill opacity=0.1, color=\inexactbfgscolor] table [x=rob-param, y=avg-krylov, col sep=comma] {data/result-rob-heart-BFGS-inexact-P2.csv};
      \legend {\newton, \inexactnewton, \inexactbfgs}
    \end{plot}
  \end{subfigure}
  \caption{\Ptwo}
  \end{subfigure}
  \caption{Robustness with respect to the data, Heartbeat test. Comparison of the nonlinear, total linear  and average linear iterations, averaged over the first 10 timesteps.}
  \label{fig:heart rob all}
\end{figure}

\subsection{Robustness with respect to the data}
To study the robustness, we varied  in each test case the parameter which resulted in a larger deformation. In the Cook test this refers to the vertical load $\tau$, in the Twist case this refers to the angle of rotation in the boundary conditions, and on the Heartbeat case this refers to the peak activation constant $C_\text{AP}$.
\begin{description}
	\item[Cook test.]
	We vary the magnitude of the distributed traction term $\tau$  on the right from $1.5\cdot10^6 \,\text{Pa}$ to $2\cdot 10^6 \,\text{Pa}$, where values over the upper bound considered make all methods diverge. The nonlinear iterations, total linear iterations and average linear iterations per nonlinear step are shown in Figure \ref{fig:cook rob all}. We have computed the solution for four values of $\tau$, and plotted the result only when the solver converged. Most notably, only the inexact Newton-Krylov method is able to solve the highest load value, confirming that it can show more robust behavior than the Newton-Krylov method. Moreover, the inexact-BFGS does not converge for some values in case of second order elements, making it the least robust method in this scenario. In general, the difference in robustness is much clearer with second order elements, as also BFGS-preonly and Newton-Krylov converge only for two of the four values considered for $\tau$. When looking at the linear iteration counts, we see that the methods that converge tend to maintain a roughly constant number of average linear iterations, so that the increase in total linear iterations is given by the additional nonlinear iterations incurred. The only exception is the Newton-Krylov method for second order elements. Interestingly, inexact-BFGS is more robust than BFGS-preonly for first order elements, but instead BFGS-preonly is slightly more robust than inexact-BFGS for second order elements. There is no advantage in considering a Newton-Krylov method, but instead its inexact variant is much more robust with respect to $\tau$.

	\item[Twist test.] We vary the angle in the boundary conditions from $\pi/6$ to $\pi/2$, and show the results in Figure \ref{fig:twist robust all}. We note that both Newton methods converge for all values under consideration, and diverge when higher ones are considered, so there is no significant difference in terms of robustness between them. Still, the Newton-Krylov method presents and increase in the total linear iterations, whereas the inexact variant shows a more robust iteration count.

    \item[Heartbeat test.] We vary the peak activation $C_{PA}$ from $10^4$ to $10^5$, and plot only the results when the methods converged. In this case, we note that Newton methods provide an improved robustness for both first and second order finite elements, and as in the other tests, the inexact variant is more robust. This effect is inverted in the BFGS methods, which are less robust and instead the BFGS-preonly shows an improved robustness with respect to the activation force. We again observe that the total linear iterations of the Newton-Krylov method increase with the problem difficulty, whereas the inexact Newton-Krylov method maintains a roughly constant number of linear iterations despite the increase of the problem difficulty.

\end{description}

\begin{table}
\begin{center}
\begin{tabular}{r|rrr|rrr|rrr|rr}
\toprule
\multicolumn{12}{c}{\Pone \ discretization, DoFs = 1576875} \\
\toprule
cores           & \multicolumn{3}{c|}{\newton}        & \multicolumn{3}{c|}{\inexactnewton}       & \multicolumn{3}{c|}{\inexactbfgs}       & \multicolumn{2}{c}{\bfgs}\\
                & nit   & lit   & $T_{sol}$       & nit   & lit   & $T_{sol}$       & nit   & lit   & $T_{sol}$       & nit   & $T_{sol}$\\
\midrule
1		& 5	& 12.4	& 718.2		& 5	& 12.0	& 708.1		& 7	& 10.4	& 516.4		& 184	& 5815.2 \\
2		& 5	& 13.0	& 616.5		& 5	& 12.0	& 505.2		& 7	& 10.7	& 330.6		& 132	& 2786.8  \\
4		& 5	& 13.2	& 297.5		& 5	& 12.6	& 293.9		& 7	& 10.1	& 173.8		& 43	& 176.2 \\
8		& 6	& 11.0	& 232.4		& 5	& 13.0	& 184.9		& 7	& 12.1	& 116.7		& 146	& 1274.1 \\
16		& 5	& 13.6	& 136.5		& 5	& 12.4	& 135.9		& 7	& 9.4	& 64.5		& 203	& 1382.1 \\
32		& 4	& 16.3	& 94.4		& 5	& 13.8	& 108.3		& 7	& 11.4	& 56.6		& 146	& 780.8  \\
\bottomrule
\toprule
\multicolumn{12}{c}{\Ptwo \ discretization, DoFs = 1576875} \\
\toprule
cores           & \multicolumn{3}{c|}{\newton}        & \multicolumn{3}{c|}{\inexactnewton}       & \multicolumn{3}{c|}{\inexactbfgs}        & \multicolumn{2}{c}{\bfgs}\\
                & nit   & lit   & $T_{sol}$     & nit   & lit   & $T_{sol}$     & nit   & lit   & $T_{sol}$     & nit   & $T_{sol}$\\
\midrule
1		& 4	& 20.3	& 1639.2	& 5	& 15.8	& 1791.4	& 7	& 13.4	& 1372.3	& 62	& 1414.9 \\
2		& 4	& 20.0	& 923.4	 	& 5	& 16.0	& 1042.3	& 7	& 13.4	& 691.0		& 62	& 721.5 \\
4		& 4	& 20.5	& 565.3		& 5	& 16.4	& 637.3		& 7	& 13.8	& 417.0		& 62	& 408.2 \\
8		& 4	& 20.8	& 363.5		& 5	& 16.8	& 409.4		& 7	& 14.6	& 263.5		& 61	& 276.6 \\
16		& 4	& 20.8	& 254.5		& 5	& 16.8	& 279.9		& 7	& 15.1	& 154.9		& 62	& 137.6 \\
32		& 4	& 21.5	& 185.9		& 5	& 17.2	& 212.2		& 7	& 14.0	& 110.9		& 63	& 126.7 \\
\bottomrule
\end{tabular}
\caption{Scalability, Cook test. nit$\,\coloneqq$ nonlinear iterations, lit$\,\coloneqq$ average linear iterations per nonlinear iteration, $T_{sol}\coloneqq$ CPU time in seconds.}
\label{tab_scal_cook}
\end{center}
\end{table}

\begin{figure}
  \begin{subfigure}{0.49\textwidth}
    \begin{loglogplot}[height=5cm]{Cores}{CPU time}
    \addplot+[fill opacity=0.1, color=\newtoncolor] table [x=cores, y=time, col sep=comma] {data/result-scalab-cook-Newton-P1.csv};
    \addplot+[fill opacity=0.1, color=\inexactnewtoncolor] table [x=cores, y=time, col sep=comma] {data/result-scalab-cook-Newton-inexact-P1.csv};
    \addplot+[fill opacity=0.1, color=\bfgscolor] table [x=cores, y=time, col sep=comma] {data/result-scalab-cook-BFGS-P1.csv};
    \addplot+[fill opacity=0.1, color=\inexactbfgscolor] table [x=cores, y=time, col sep=comma] {data/result-scalab-cook-BFGS-inexact-P1.csv};
    \legend{\newton, \inexactnewton, \bfgs, \inexactbfgs}
  \end{loglogplot}
    \caption{\Pone}
  \end{subfigure}
  \begin{subfigure}{0.49\textwidth}
    \begin{loglogplot}[height=5cm]{Cores}{}
    \addplot+[fill opacity=0.1, color=\newtoncolor] table [x=cores, y=time, col sep=comma] {data/result-scalab-cook-Newton-P2.csv};
    \addplot+[fill opacity=0.1, color=\inexactnewtoncolor] table [x=cores, y=time, col sep=comma] {data/result-scalab-cook-Newton-inexact-P2.csv};
    \addplot+[fill opacity=0.1, color=\bfgscolor] table [x=cores, y=time, col sep=comma] {data/result-scalab-cook-BFGS-P2.csv};
    \addplot+[fill opacity=0.1, color=\inexactbfgscolor] table [x=cores, y=time, col sep=comma] {data/result-scalab-cook-BFGS-inexact-P2.csv};
    \legend{\newton, \inexactnewton, \bfgs, \inexactbfgs}
  \end{loglogplot}
    \caption{\Ptwo}
  \end{subfigure}
  \caption{Scalability, Cook test. CPU times with respect the number of cores. }
  \label{fig:cook scalab time}
\end{figure}

\begin{table}
\begin{center}
\begin{tabular}{r|rrr|rrr}
\toprule
\multicolumn{7}{c}{DoFs = 519747} \\
\toprule
cores           & \multicolumn{3}{c}{\newton}         & \multicolumn{3}{c}{\inexactnewton} \\
                & nit   & lit   & $T_{sol}$       & nit   & lit   & $T_{sol}$ \\
\midrule
1		& 6	& 800.0	& $3.68\cdot10^4$	& 16	& 240.4	& $2.61\cdot10^4$ \\
2		& 6	& 806.2	& $2.13\cdot10^4$ 	& 16	& 225.3	& $1.62\cdot10^4$ \\
4		& 6	& 806.0	& $1.37\cdot10^4$ 	& 16	& 222.4	& $1.06\cdot10^4$ \\
8		& 6	& 826.3	& $8.02\cdot10^3$	& 16	& 237.3	& $6.56\cdot10^3$ \\
16		& 6	& 835.5	& $4.62\cdot10^3$	& 16	& 240.3	& $3.94\cdot10^3$ \\
32		& 6	& 828.5	& $3.44\cdot10^3$	& 16	& 244.8	& $3.10\cdot10^3$ \\
\bottomrule
\end{tabular}
\caption{Scalability, Twist test. nit$\,\coloneqq$ nonlinear iterations, lit$\,\coloneqq$ average linear iterations per nonlinear iteration, $T_{sol}\coloneqq$ CPU time in seconds.}
\label{tab_scal_twist}
\end{center}
\end{table}

\begin{figure}
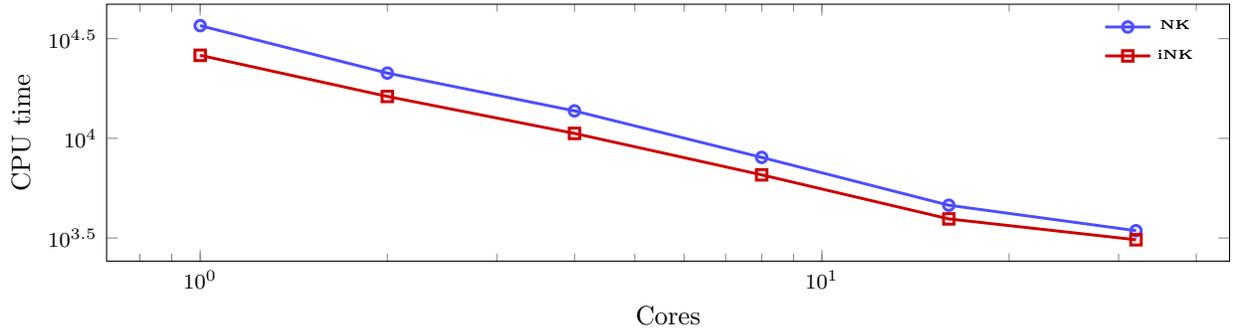

      \begin{loglogplot}[height=5cm]{Cores}{CPU time}
      \addplot+[fill opacity=0.1, color=\newtoncolor] table [x=cores, y=time, col sep=comma] {data/result-scalab-twist-Newton-P2.csv};
      \addplot+[fill opacity=0.1, color=\inexactnewtoncolor] table [x=cores, y=time, col sep=comma] {data/result-scalab-twist-Newton-inexact-P2.csv};
      \legend {\newton, \inexactnewton}
    \end{loglogplot}
  \caption{Scalability, Twist test. CPU time with respect to the number of cores.}
  \label{fig:twist scalab}
\end{figure}

\begin{table}
\begin{center}
\begin{tabular}{r|rrr|rrr|rrr|rr}
\toprule
\multicolumn{12}{c}{\Pone \ discretization, DoFs = 149511} \\
\toprule
cores           & \multicolumn{3}{c|}{\newton}        & \multicolumn{3}{c|}{\inexactnewton}       & \multicolumn{3}{c|}{\inexactbfgs}       & \multicolumn{2}{c}{\bfgs}\\
                & nit   & lit   & $T_{sol}$       & nit   & lit   & $T_{sol}$       & nit   & lit   & $T_{sol}$       & nit   & $T_{sol}$\\
\midrule
1		& 4.5	& 9.6	& 111.0		& 5.0	& 6.0	& 113.6		& 10.9	& 2.4	& 69.4		& 23.2	& 87.7 \\
2		& 4.5	& 9.6	& 60.9		& 5.0	& 6.0	& 62.3		& 11.0	& 2.5	& 30.1		& 23.2	& 45.5 \\
4		& 4.5	& 9.6	& 26.5		& 5.0	& 6.0	& 33.2		& 10.8	& 2.7	& 18.8		& 23.3	& 23.3 \\
8		& 4.5	& 9.8	& 18.2		& 5.0	& 6.2	& 18.9		& 11.0	& 2.6	& 10.3		& 23.9	& 12.6 \\
16		& 4.5	& 9.7	& 11.1		& 5.0	& 6.1	& 11.5		& 10.8	& 2.6	& 4.8		& 23.5	& 6.8 \\
32		& 4.5	& 9.7	& 7.6		& 5.0	& 6.2	& 7.9		& 10.9	& 2.8	& 3.1		& 23.6	& 4.0 \\
\bottomrule
\toprule
\multicolumn{12}{c}{\Ptwo \ discretization, DoFs = 1123515} \\
\toprule
cores           & \multicolumn{3}{c|}{\newton}        & \multicolumn{3}{c|}{\inexactnewton}       & \multicolumn{3}{c|}{\inexactbfgs}        & \multicolumn{2}{c}{\bfgs}\\
                & nit   & lit   & $T_{sol}$     & nit   & lit   & $T_{sol}$     & nit   & lit   & $T_{sol}$     & nit   & $T_{sol}$\\
\midrule
1		& 4.6	& 13.4	& 1061.0	& 5.3	& 6.5	& 965.4		& 14.1	& 3.5	& 637.4		& 30.2	& 515.1 \\
2		& 4.6	& 13.5	& 512.8		& 5.3	& 6.4	& 569.0		& 13.4	& 3.6	& 287.4		& 30.4	& 287.7 \\
4		& 4.6	& 13.5	& 324.5		& 5.0	& 6.5	& 294.3		& 15.3	& 3.4	& 184.9		& 30.5	& 149.8 \\
8		& 4.6	& 13.4	& 184.1		& 5.3	& 6.6	& 176.0		& 14.5	& 3.6	& 100.7		& 30.1	& 81.7 \\
16		& 4.6	& 13.4	& 113.2		& 5.2	& 6.6	& 107.8		& 14.5	& 3.3	& 59.1		& 30.5	& 40.5 \\
32		& 4.6	& 13.4	& 82.0		& 5.1	& 6.5	& 74.5 		& 14.2	& 3.5	& 43.3		& 30.3	& 33.6 \\
\bottomrule
\end{tabular}
\caption{Scalability, Heartbeat test. nit$\,\coloneqq$ nonlinear iterations, lit$\,\coloneqq$ average linear iterations per nonlinear iteration, $T_{sol}\coloneqq$ CPU time in seconds.}
\label{tab_scal_heart}
\end{center}
\end{table}

\begin{figure}
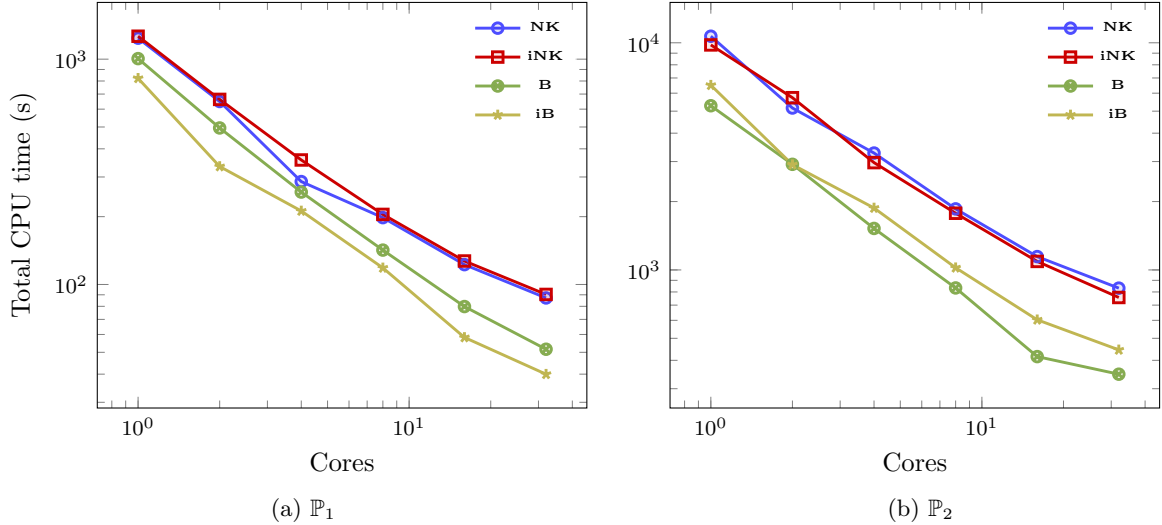

  \begin{subfigure}{0.49\textwidth}
    \begin{loglogplot}{Cores}{Total CPU time (s)}
    \addplot+[fill opacity=0.1, color=\newtoncolor] coordinates {
    (1,1239.2973704859614)
    (2,648.8568870108575)
    (4,286.2344272714108)
    (8,198.02774427458644)
    (16,122.48819336481392)
    (32,87.15087462030351)
    };
    \addplot+[fill opacity=0.1, color=\inexactnewtoncolor] coordinates {
    (1,1264.1495023109019)
    (2,662.972876381129)
    (4,356.9405936654657)
    (8,204.6662153713405)
    (16,127.11849938705564)
    (32,90.33599272370338)
    };
    \addplot+[fill opacity=0.1, color=\bfgscolor] coordinates {
    (1,1005.6209271810949)
    (2,495.1001786123961)
    (4,257.13379335962236)
    (8,142.11290678568184)
    (16,79.81870525889099)
    (32,51.59551910869777)
    };
    \addplot+[fill opacity=0.1, color=\inexactbfgscolor] coordinates {
    (1,822.2470116131008)
    (2,333.9728074353188)
    (4,212.00657227635384)
    (8,118.49205537326634)
    (16,58.209598591551185)
    (32,39.983160100877285)
    };
    \legend{\newton, \inexactnewton, \bfgs, \inexactbfgs}
  \end{loglogplot}
    \caption{\Pone}
  \end{subfigure}
  \begin{subfigure}{0.49\textwidth}
    \begin{loglogplot}{Cores}{}
    \addplot+[fill opacity=0.1, color=\newtoncolor] coordinates {
    (1,10687.74784265086)
    (2,5162.593246784061)
    (4,3270.4420484751463)
    (8,1857.0429484266788)
    (16,1144.8118587359786)
    (32,831.467165466398)
    };
    \addplot+[fill opacity=0.1, color=\inexactnewtoncolor] coordinates {
    (1,9786.219681892544)
    (2,5731.672663629055)
    (4,2968.2274741865695)
    (8,1776.6269434746355)
    (16,1090.3236375711858)
    (32,756.469524603337)
    };
    \addplot+[fill opacity=0.1, color=\bfgscolor] coordinates {
    (1,5283.185650136322)
    (2,2917.24963394925)
    (4,1523.0937974024564)
    (8,833.1856028418988)
    (16,415.40010945498943)
    (32,347.5394565425813)
    };
    \addplot+[fill opacity=0.1, color=\inexactbfgscolor] coordinates {
    (1,6505.937527783215)
    (2,2907.9546362850815)
    (4,1874.026005603373)
    (8,1022.5915530901402)
    (16,603.3599721081555)
    (32,444.7639068197459)
    };
    \legend{\newton, \inexactnewton, \bfgs, \inexactbfgs}
  \end{loglogplot}
    \caption{\Ptwo}
  \end{subfigure}

  \caption{Scalability, Heartbeat test. CPU times with respect to the number of cores. }
  \label{fig:heart scalab strong}
\end{figure}

\subsection{Scalability}
We test the strong scalability with respect to the finest mesh used in the performance tests. In the heartbeat test, we consider only the first 10 timesteps. 
\begin{description}
	\item[Cook test.] We report the nonlinear and average linear iterations and the CPU times in Table \ref{tab_scal_cook}. The CPU times are also displayed in Figure \ref{fig:cook scalab time}. We note that all methods are robust with respect to the number of cores, as could have been expected in virtue of the robustness of the preconditioner. Still, it is very interesting to see that there is a point in which BFGS-preonly converges with a much lower iteration number in the \Pone case, namely when using 8 processors. This is the only case in which the CPU time is comparable to that of the inexact-BFGS method. For any other number of cores, the inexact-BFGS outperforms all other methods for both \Pone and \Ptwo \,elements. As already observed in the performance test, the BFGS-preonly method is less robust than the other methods.

	\item[Twist test.] We report the nonlinear and average linear iteration counts, together with the CPU times in Table \ref{tab_scal_twist}. The CPU times are also displayed in Figure \ref{fig:twist scalab}. We note that this test in particular is in complete agreement with the previous ones, meaning that both Newton-Krylov methods perform similarly and exhibit an overall robust behavior. Both methods present adequate strong scaling, with the inexact method showing better solution times.

	\item[Heartbeat test.] For this test, we report the nonlinear and average linear iterations, and CPU times averaged for the first 10 timesteps in Table \ref{tab_scal_heart}. The total CPU times as function of the number of cores are displayed in Figure \ref{fig:heart scalab strong}. We can appreciate that all methods are robust with respect to the number of cores, as these yield no significant variations. From the scalability curves shown in Figure \ref{fig:heart scalab strong}, we note that all methods show adequate scaling, with the BFGS methods presenting an overall better performance. More specifically, inexact-BFGS is faster for \Pone, and BFGS-preonly is faster for \Ptwo.
\end{description}

\section{Conclusions}\label{section:discussion}
This work presents a detailed numerical study of some of the main numerical solvers used in large scale simulations of nonlinear mechanics, where all of the fundamental components of the methods are included: the nonlinear solver, the linear solver and the preconditioner. The main finding from this study is that superlinear solvers such as inexact Newton-Krylov and BFGS methods present an overall  robustness similar to the gold standard  represented by the Newton-Krylov method in nonlinear elasticity, but with improved solution times. 

In difficult problems such as static and incompressible mechanics, the choice of the solver can be defined a-priori according to the convergence hypotheses required for each method, where in fact we have seen the BFGS methods fail for incompressible mechanics, as they usually require some form of convexity to guarantee convergence. Indeed, BFGS methods showed the best solution times in static mechanics, while inexact Newton-Krylov showed the best performance for the incompressible tests. The time-dependent heartbeat test showed us that the convexity contributed by the inertia term can further improve the performance of the BFGS method, which showed much more robust iteration counts in this scenario than those shown in the static mechanics case. This comes as a surprise, given that the constitutive model used in the heartbeat test is much more complex than the models of the other tests. These nonlinear solvers present no drawbacks in terms of both performance and strong scalability, suggesting that in the context of nonlinear elasticity inexact Newton-Krylov and quasi-Newton methods should be preferred to standard Newton-Krylov methods.

We conclude by mentioning two possible extensions of this work. The first one is the inclusion in this study of Jacobian lagging techniques \cite{brown2013low}, which consist in avoiding the reassembly of the Jacobian for some Newton iterations. These lagging techniques could potentially further improve the performance of the Newton methods considered. Instead, inspired by our inexactness taxonomy, we note that a Jacobian reassembly could be considered for the BFGS methods, where the Jacobian gets sometimes reassembled to improve the quality of the initial Jacobian. Both strategies are non trivial to implement, and require a detailed problem-specific study to obtain the best performance. The second interesting extension is the inclusion in our nonlinear mechanics study of a polyconvex potential \cite{ball1976convexity}, even if solvers for this specific type of potential are not yet available, to the best of our knowledge. One interesting work in this direction is \cite{bonet2015computational}, where a Hu-Washizu \cite{washizu1968variational} formulation is used to exploit the convexity of the auxiliary function arising from the polyconvex potential. Still, this formulation results in a problem with many auxiliary variables, and its efficient numerical approximation remains a research topic.

\appendix
\section{PETSc options}\label{appendix:petsc}
The commands used to invoke each of the methods described in Section \ref{section:solvers} is detailed in what follows.

\begin{itemize}
    \item {\bf Newton-Krylov } (\newton)
    \begin{lstlisting}[language=bash, frame=single, caption=PETSc commands to use Newton-Krylov.]
    -snes_type newtonls
    -ksp_type  gmres  # minres should work as well
    -pc_type   hypre
    -ksp_atol  1e-10
    -ksp_rtol  1e-6
    -snes_linesearch_type basic
    \end{lstlisting}
    \item {\bf Inexact Newton-Krylov} (\inexactnewton)
    \begin{lstlisting}[language=bash, frame=single, caption=PETSc commands to use inexact Newton-Krylov.]
    -snes_type newtonls
    -ksp_type  gmres
    -pc_type   hypre
    -ksp_atol  1e-14
    -ksp_rtol  1e-2  # Should not be used, set for safety
    -snes_ksp_ew
    -snes_ksp_ew_rtol0   1e-1
    -snes_ksp_ew_rtolmax 0.1
    -snes_linesearch_type basic
    \end{lstlisting}
    \item {\bf Inexact-BFGS} (\inexactbfgs)
    \begin{lstlisting}[language=bash, frame=single, caption=PETSc commands to use inexact-BFGS.]
    -snes_type qn
    -ksp_type  gmres 
    -pc_type   hypre
    -ksp_atol  1e-14
    -ksp_rtol  1e-2
    -snes_qn_type lbfgs
    -snes_qn_m 50 
    -snes_qn_scale_type jacobian
    -snes_linesearch_type basic
    \end{lstlisting}
    \item {\bf BFGS-preonly} (\bfgs)
    \begin{lstlisting}[language=bash, frame=single, caption=PETSc commands to use BFGS.]
    -snes_type qn
    -ksp_type  preonly
    -pc_type   hypre
    -snes_qn_type lbfgs
    -snes_qn_m 50
    -snes_qn_scale_type jacobian
    -snes_linesearch_type basic
    \end{lstlisting}
\end{itemize}

\section{Schur complement preconditioners}\label{appendix:schur}
A Schur complement preconditioner is one that arises from a block LU factorization according to two (arbitrary) index sets of a matrix. If we consider a block matrix $\mathbf M$ given by the general structure
    \[ \mathbf M = \begin{bmatrix} \mathbf A & \mathbf B_1 \\ \mathbf B_2 & \mathbf C \end{bmatrix} ,\]
with $\mathbf A$ invertible, a Gaussian elimination procedure yields
    \begin{equation}\label{eq:schur}
    \mathbf M = \begin{bmatrix} I & \mathbf 0 \\ \mathbf B_2\mathbf A^{-1} \end{bmatrix} \begin{bmatrix} \mathbf A & \mathbf 0 \\ \mathbf 0 & \mathbf C -\mathbf B_2\mathbf A^{-1}\mathbf B_1  \end{bmatrix} \begin{bmatrix} \mathbf I & \mathbf A^{-1}\mathbf B_1 \\ \mathbf 0 & \mathbf I \end{bmatrix}.
         \end{equation}
We note that if $\mathbf C$ is invertible, the same procedure can be applied with respect to it. Schur complement based preconditioners enjoy excellent theoretical properties, as the preconditioned system possesses at most 3 distinct eigenvalues \cite{murphy2000note}, implying that it converges in at most 3 iterations of a Krylov subspace method. This is true whenever the Schur complement $\mathbf S = \mathbf C -\mathbf B_2\mathbf A^{-1}\mathbf B_1$ is evaluated exactly, which is usually computationally intractable. One simple approximation, which is the one we use for the preconditioning the incompressible mechanics problem, is to consider the approximation $\mathbf A^{-1} \approx \text{diag}^{-1}\,(\mathbf A)$. Another important point is that in general using all three blocks arising from the factorization in \eqref{eq:schur} is not necessary, and instead it is sufficient to consider a lower factorization (first two blocks) or an upper factorization (last two blocks). Whenever $\mathbf C=\mathbf 0$, it is also possible to use a diagonal factorization (middle block only).

\section*{Acknowledgments}
N. Barnafi and L. F. Pavarino have been supported by
grants of MIUR (PRIN 2017AXL54F$\_$002) and INdAM--GNCS. 
N. Barnafi and S. Scacchi have been supported by grants of MIUR (PRIN
2017AXL54F$\_$003) and INdAM-GNCS. 
The Authors are also grateful to the University of Pavia, the University of Milan, and the CINECA laboratory for the usage of the EOS, INDACO and Galileo100 clusters, respectively.

\bibliography{main}
\bibliographystyle{alpha}

\end{document}